\input amstex
\documentstyle{amsppt}

\bigskip

\centerline{\bf MOON-TYPE THEOREMS ON CIRCUITS IN STRONGLY}

\smallskip

\centerline{\bf CONNECTED TOURNAMENTS  OF ORDER $\bold{N}$ AND DIAMETER $\bold{D}$}

\bigskip

\centerline{\bf S.V. Savchenko}

\bigskip

\centerline{ L.D. Landau Institute for Theoretical Physics, Russian Academy of Sciences}
\centerline{ Kosygin str. 2, Moscow 119334, Russia}

\smallskip

\centerline{\sl In memory of Valerii Evgen$'$evich Tarakanov}

\bigskip

Let $T$ be a strongly connected tournament of order $n\ge 4$ whose
diameter does not exceed $d\ge 3.$ Denote by $c_{\ell}(T)$
the number of circuits of length $\ell$ in $T.$
In our recent paper, we construct a strongly connected tournament  $T_{d,n}$
of order $n$ with diameter $d$ and conjecture that
$c_{\ell}(T)\ge c_{\ell}(T_{d,n})$ for any $\ell=3,...,n.$
In particular, for $d=n-1,$
this inequality is true and yields the known Moon (lower) bound
$c_{\ell}(T)\ge n-\ell+1.$
Moreover, we suggest that
if $n+3\le 2d,$ then for any given $\ell$ taken in the range $n-d+3,...,d,$
the equality $c_{\ell}(T)=c_{\ell}(T_{d,n})$
implies that $T$ is isomorphic to
$T_{d,n}$ or its converse $T_{d,n}^{-}.$
For $d=n-1,$ the corresponding particular statement
is nothing else than
Las Vergnas' theorem. Recently,
we have confirmed the posed conjecture for the case $d=n-2.$
In the present paper, we show that it is also true
for $d=n-3.$

\smallskip

Bibliography: 18 titles.

\smallskip


\smallskip


\bigskip

\centerline{\bf 1. Introduction}

\smallskip

By definition, a {\sl tournament} $T$ of order $n$
is an orientation of
the complete graph $K_{n}$ on $n$ vertices.
In other words, any two vertices $v$ and $w$ are joined by exactly one of the possible arcs
$(v,w)$ and $(w,v)$.
For describing $T,$
we shall use a standard terminology
(see, for instance, [1]).
Let $V(T)$ be its vertex-set and $A(T)$ be its arc-set.
If $(v,w)\in A(T),$ then
we shall say that the vertex $v$ {\sl dominates} the vertex $w$
in $T,$ and write $v\to w.$
In turn, for two tournaments
$T_{i}$ and $T_{j}$ whose vertex-sets $V(T_{i})$ and $V(T_{j})$
are disjoint, the notation
$T_{i}\Rightarrow T_{j}$ means that each vertex in $T_{i}$ dominates
each vertex in $T_{j}.$
For $n$ vertex-disjoint
tournaments $T_{1},...,T_{n},$
the {\sl composition } $T(T_{1},...,T_{n})$
is the tournament obtained from $T$ by replacing its vertices $w_{1},...,w_{n}$
and the binary relation $\to$ between them
with the tournaments $T_{1},...,T_{n}$ and the binary relation
$\Rightarrow$ between the last, respectively.
If $w_{i}\in V(T_{i}),$ where $i=1,...,n,$ we will say that
the remaining vertices of $T_{i}$ {\sl enlarges} $w_{i}$
in $T(T_{1},...,T_{n}).$

Without any doubt, tournaments form the best studied class
of directed graphs. For the convenience of the reader,
the known results on tournaments
are indexed by the Roman numerals
and the theorems proved by the author himself
are numbered with the Arabic numerals,
as usual. The first classical result
concerns {\sl hamiltonian paths} which, by definition, contain
all vertices of $T.$

\smallskip

{\bf Theorem I [13].} {\sl Any tournament
contains an odd number of hamiltonian paths.}

\smallskip

In particular, every tournament admits at least
one hamiltonian path. This lower bound is sharp.
For any order $n,$ it is attained only at the {\sl transitive}
tournament $TT_{n}$ of order $n.$
By definition, if three vertices $v,u,$ and $w$ of it are related
by $v\to u\to w,$ then $v\to w.$
So, if $z_{0},...,z_{n-1}$ is its hamiltonian path,
then $z_{i}\to z_{j}$ for $i<j.$ We denote this particular (labeled)
tournament by $TT_{n}(z_{0},...,z_{n-1}).$
Obviously, $TT_{n}$
contains no circuits at all. By this reason, it is also often
called {\sl acyclic}. (Note that this property as well as
the transitivity condition completely characterizes
$TT_{n}$ among all tournaments of order $n$.)

Let $N^{-}_{T}(w)$ be the set of vertices dominating $w$ in $T$
and $N^{+}_{T}(w)$ be the set of vertices dominated by $w$ in $T.$
The quantities $d^{+}_{T}(w)=|N^{+}_{T}(w)|$ and $d^{-}_{T}(w)=|N^{-}_{T}(w)|$
are called
the {\sl out-degree} and {\sl in-degree} of the vertex $w,$ respectively.
Obviously, $d^{+}_{T}(w)=|T|-1-d^{-}_{T}(w).$
Hence, in the sequel, we mainly deal with the out-degrees of vertices.
For instance, if $T=TT_{n}(z_{0},...,z_{n-1}),$ then
$d^{+}_{T}(z_{i})=n-1-i$ for $i=0,...,n-1.$
Note that the sequence $0,...,n-1$ of the out-degrees of vertices
uniquely determines $TT_{n}$ in the class of all tournaments of order $n.$

In the sequel, we mainly consider
{\sl strongly connected} (or, merely, {\sl strong}) tournaments $T.$
This means that for any two vertices $v$ and $w$ of $T,$
there exists a path from
$v$ to $w$ in $T.$ Camion's theorem [5] states that
any strongly connected tournament admits a {\sl hamiltonian circuit}.
In fact, the following lower bound on the number $c_{\ell}(T)$
of circuits of length $\ell$ in such a $T$ of order $n$ even holds.

\smallskip

{\bf Theorem II [8].} {\sl Let $T$ be a strongly connected
tournament of order $n\ge 3.$
Then for $3\le \ell\le n,$ the inequality
$c_{\ell}(T)\ge n-\ell+1$ holds.}

\smallskip

This lower bound
is also sharp. In particular, it is
attained at the strongly connected tournament $T_{n-1,n}(z_{0},...,z_{n-1})$
with hamiltonian path $z_{0},...,z_{n-1}$ such that
$z_{j}\to z_{i}$ for $j>i+1.$
If $n=1,2,$ then we can assume that
$T_{n-1,n}=TT_{n}.$ For $n=3,$ the tournament
$T_{n-1,n}$ is the {\sl cyclic triple} $\Delta$
and for $n=4,$ it is obtained from $\Delta$
by replacing one of its vertices with $TT_{2}.$
Note that if $T=T_{n-1,n}(z_{0},...,z_{n-1}),$ then
$d_{T}^{+}(z_{0})=1,$ $d_{T}^{+}(z_{i})=i$ for $i=1,...,n-2,$
and $d_{T}^{+}(z_{n-1})=n-2.$
However, the sequence $1,1,2,...,n-3,n-2,n-2$ of the out-degrees of vertices
does not uniquely determine $T_{n-1,n}$
in the class of all tournaments of order $n.$
As it was shown in [3], for $\ell=3,$
the lower bound of Theorem II is attained only at strong tournaments
with this out-degree sequence.
They all also contain exactly one hamiltonian circuit.
The class of tournaments with this property
(it is a bit wider than the family consisting of all strongly connected
tournaments with the least possible number of cycles of length $3$)
was completely described in [6].
However, according to the Las Vergnas theorem [7],
for $\ell=4,...,n-1,$  there are no other strongly
connected tournaments of order $n$
but $T_{n-1,n}$ for which the lower bound of Theorem II
on $c_{\ell}(T)$ is achieved.

\smallskip

{\bf Theorem III [7].} {\sl If $T$ is a strongly connected tournament
of order $n\ge 5$ and $c_{\ell}(T)=n-\ell+1$ for some $\ell=4,...,n-1,$
then $T$ is isomorphic to $T_{n-1,n}.$}

\smallskip

Note that all the above classical results
on a strong tournament $T$ involve only one given
parameter of it, namely, the order $n.$
In our study started in [16], we also restrict its diameter.
By definition, the {\sl distance} $d_{T}(x,y)$
between two vertices $x$ and $y$ in $T$
is the length
of a shortest path from $x$ to $y$ in $T.$
In turn, the {\sl diameter} of the tournament $T$ is
the maximum possible distance between two distinct vertices in $T.$
Obviously, the diameter cannot be greater than $n-1$
and the digraph $T_{n-1,n}$ introduced above
is the unique strongly connected
tournament of order $n$ and diameter $n-1.$

Now, let ${\Cal T}_{d,n}$ be the class
of all strongly connected tournaments
of order $n$ and diameter $d.$
For $n>d\ge 3,$ in the class ${\Cal T}_{d,n},$
we select two tournaments
$$T_{d,n}=T_{d,d+1}(TT_{\lfloor\frac{n-d+1}{2}\rfloor},v_{1},...,v_{d-1},
TT_{\lceil\frac{n-d+1}{2}\rceil})$$
and
$$T_{d,n}^{-}=T_{d,d+1}(TT_{\lceil\frac{n-d+1}{2}\rceil},v_{1},...,v_{d-1},
TT_{\lfloor\frac{n-d+1}{2}\rfloor}),$$
where $\lceil \frac{n-d+1}{2}  \rceil$ is the smallest integer greater than or equal to $\frac{n-d+1}{2}$
and $\lfloor \frac{n-d+1}{2}\rfloor$ is the greatest integer smaller than or equal to $\frac{n-d+1}{2}.$

In the sequel, we say that a tournament
is the {\sl converse} of another one if the former is obtained
by reversing all arcs in the latter.
Obviously, $T_{d,n}^{-}$ is isomorphic to the converse of $T_{d,n}.$
If $n-d+1$ is even, then $T_{d,n}^{-}$ coincides with
$T_{d,n}.$
Note also that for any (not only even) $n\ge 4,$ we have
$T_{3,n}\cong T_{3,n}^{-}\cong \Delta(v_{1},v_{2},TT_{n-2}),$
where $\Delta$ is the cyclic triple with the hamiltonian circuit
$v_{1},v_{2},v_{3},v_{1}.$
In the remaining cases, $T_{d,n}^{-}$ is not isomorphic to $T_{d,n}.$

\smallskip

{\bf Conjecture 1 [16].} {\sl
Let $T$ be a strongly connected tournament of order $n\ge 4$
whose diameter does not exceed $d\ge 3$.
Then
$c_{\ell}(T)\ge c_{\ell}(T_{d,n})$
for each $\ell=3,...,n.$
Moreover, if $d\ge \frac{n+3}{2},$
then for given $\ell$ in the range $n-d+3,...,d,$
the equality $c_{\ell}(T)=c_{\ell}(T_{d,n})$
implies that
$T$ is isomorphic to
$T_{d,n}$ or to its converse $T_{d,n}^{-}.$}

\smallskip

Note that for $2d\ge n+3$ and $n-d+3\le \ell \le d,$
we have
$$c_{\ell}(T_{d,n})=d-\ell+2^{\lceil\frac{n-d+1}{2}\rceil}+
2^{\lfloor\frac{n-d+1}{2}\rfloor}-2.$$
However,
if $n+d-1\le 2\ell$ (in particular, this means that $d\le \ell$),
then
$$c_{\ell}(T_{d,n})=\binom{n-d+1}{n-\ell}.$$
The reader can understand how to deduce
the expressions for $c_{\ell}(T_{d,n})$ in the cases
$n+d-1\le 2\ell$ and $n-d+3\le \ell \le d$
if he reads Remark 4 and
the text in the lines between the statements of
Theorem 1 and Conjecture 2, respectively.
In fact, the number $c_{\ell}(T_{d,n})$ can be determined for all
possible values of $\ell,d,$ and $n.$ Unfortunately, the corresponding
expression is too large for convenient presentation here.
So, we omit it.

Obviously, for $d=n-1,$  the last statement of Conjecture 1
is nothing but Las Vergnas' theorem mentioned above as Theorem III.
The case $d=n-2$ was considered in [16].
In the present paper, we also confirm the conjecture for $d=n-3.$
The corresponding proof is given in Section 3.
But before giving it, we present
some necessary results on the number $c_{\ell}(T,w)$
of cycles of length $\ell$ including a vertex $w$ in $T.$

\bigskip

\centerline{\bf 2. Moon's vertex-pancyclic theorem and non-critical vertices }

\smallskip

Standard proof of Theorems II and III
is given by induction on $n$ and uses
the existence of a vertex $w$ such that the one-vertex-deleted
subtournament $T-w$
of order $n-1$ is also strongly connected.
Such a vertex is called {\sl non-critical}.
In the opposite case (i.e. $T-w$ is not strongly connected),
it is a {\sl critical } vertex.

Obviously, to use the induction on $n,$
we need to have $c_{\ell}(T,w)\ge 1$ for some non-critical vertex $w.$
Fortunately, this is already guaranteed by Moon's vertex-pancyclic theorem
which is, without any doubt,
one of the main results in tournament theory.

\smallskip

{\bf Theorem IV [8], [9].} {\sl For any $\ell=3,...,n,$ each vertex of
a strongly connected tournament
of order $n\ge 3$ is contained in a circuit of length $\ell.$}

\smallskip

In turn, a particular statement of Theorem IV for the case $\ell=n-1$
implies the following proposition on non-critical vertices.

\smallskip

{\bf Corollary I [8].} {\sl Any strongly connected tournament $T$
of order $n\ge 4$ contains at least
two non-critical vertices.}

\smallskip

{\bf Remark 1}. Frankly speaking,
the statement of Corollary I can be proved directly,
without the use of Theorem IV (for instance, see [12]). Moreover, it can be extended to
the class of oriented digraphs, to which tournaments also belong. Recall that
a digraph $D$ is {\sl oriented} if it admits no cycles of length $2.$
By definition,
the degree $d_{D}(x)$ of a vertex $x$ is the number of vertices $y$ such that
either $y\to x$ or $x\to y.$ Obviously, $d_{T}(x)=n-1$ for each vertex $x$
in any tournament
$T$ of order $n.$ In [14], the author proved that if the minimum
degree is not less than $\frac{3}{4}n,$ where $n\ge 4,$
then a strongly connected oriented digraph $D$ of order $n$ admits at least
two non-critical vertices. However, it was shown recently in [11] that
the same holds if the weaker
condition $d_{D}(x)\ge \frac{n+2}{2}$ is satisfied for each vertex $x$
in $D$. The example of a standard orientation of a complete bipartite
graph with equal parts
implies that the coefficient of $n$ in the condition on $d_{D}(x)$
cannot be reduced and hence,
the result of [11] is best possible.

\smallskip

The lower bound of Corollary I is also sharp. Indeed, the
strongly connected tournament $T_{n-1,n}(z_{0},...,z_{n-1}),$
where $n\ge 4,$
contains exactly two non-critical vertices, namely, $z_{0}$ and $z_{n-1}.$
(Obviously, $T_{n-1,n}-z_{0}\cong T_{n-1,n}-z_{n-1}\cong T_{n-2,n-1}$.)
In the sequel, we call $z_{0}$ and $z_{n-1}$ the left and right
non-critical vertices of $T_{n-1,n},$ respectively.
It turns out that an extremal tournament is also unique by virtue of the following theorem by Thomassen.

\smallskip

{\bf Theorem V [17].} {\sl A strongly connected
tournament $T$ of order $n\ge 4$
contains exactly two non-critical vertices if and only if
$T$ is isomorphic to $T_{n-1,n}.$}

\smallskip

Let $T_{ncr}$ be the subtournament induced by the
set of non-critical vertices of $T.$
Corollary I and
Theorem V are simple consequences of the following lemma on $T_{ncr}$
which will be also useful in the sequel.

\smallskip

{\bf Lemma 1 [16].}
{\sl Let $T$ be a strongly connected tournament of order $n\ge 4$.
If $T_{ncr}$ is also strongly connected,
then either $T=T_{ncr}$ or
$T=\Delta(v_{1},v_{2},T_{ncr}).$
Otherwise, $T_{ncr}$
is not contained in any proper strongly connected subtournament of $T$.}

\smallskip

The fact that $v_{1}$ and $v_{2}$ are critical vertices in
$\Delta(v_{1},v_{2},T_{ncr}),$ indeed also follows from the known general
results on non-critical vertices in compositions of tournaments.

\smallskip
{\bf Lemma 2 [2], [15].} {\sl
For $n>1,$ the composition $T(T_{1},...,T_{n})$
is strongly connected iff the original tournament $T$ of order $n$ is so.
For this case, the set of non-critical
vertices in $T(T_{1},...,T_{n})$ consists of the elements  of
the vertex-sets of
those $T_{r}$ for which either $|T_{r}|\ge 2$ or $|T_{r}|=1$
and the corresponding vertex $w_{r}$ of the tournament $T$
is contained in $T_{ncr}.$}

\smallskip

It is not difficult to check that $c_{\ell}(T_{n-1,n},w)=1$
for each possible length $\ell\ge 3$ and any non-critical vertex $w$ of
the tournament $T_{n-1,n}.$
Strange as it may seem, this property is also characteristic for $T_{n-1,n}$
if $4\le \ell \le n-1,$ where $n\ge 5$
(note that in the case of $n=4,$
this fact is trivial for each of the two possible values
$\ell=3$ and $\ell=4$
because $T_{3,4}$ is the unique strongly connected tournament of order $4$).

\smallskip

{\bf Theorem 1.} {\sl Let $T$ be a strongly connected
tournament of order $n\ge 5$
and $4\le \ell \le n-1.$
Assume that $c_{\ell}(T,w)=1$ for any non-critical vertex $w$ of $T.$
Then $T\cong T_{n-1,n}$.}

\smallskip

{\sl Proof.}
Suppose that $c_{\ell}(T,w)=1$ for any non-critical vertex $w$ in $T.$
Note at once that if a circuit of length $\ell$
contains two non-critical vertices $w_{1}$ and $w_{2},$
then $c_{\ell}(T,w_{2})\ge 2$
because,
by the Moon vertex-pancyclic theorem,
the strongly connected subtournament $T-w_{1}$
also has a circuit of length $\ell$ including the vertex $w_{2}.$
In particular, this means that we can assume that $|T_{cr}|\ge \ell-1,$
where $T_{cr}$ is the subtournament induced by the critical vertices of $T.$
If $\ell\ge 4,$ then $|T_{cr}|\ge 3.$
Hence, Lemma 1 implies that
$T_{ncr}$ is not contained in any proper strongly connected subtournament of $T$.

By Theorem I,  the subtournament $T_{ncr}$ admits a hamiltonian
path $P=w_{1},...,$ $w_{p}$. Note that since $|T_{ncr}|\ge 2,$
we have $p\ge 2$ and hence, $w_{p}\neq w_{1}.$
Let $Q$ be a shortest path from $w_{p}$ to $w_{1}$ in $T$
and $T(Q)$ be the subtournament induced by the vertex-set of $Q.$
The fact that $Q$ is a shortest path in $T$ and $|Q|\ge 3$ means that
$T(Q)$ is strongly connected and
the diameter of $T(Q)$ equals $|Q|-1.$ As we have already
noticed in the introduction, this condition uniquely determines
a strongly connected tournament: $T(Q)\cong T_{|Q|-1,|Q|}.$
The concatenation of $P$ and $Q$ forms a closed
walk $W.$ Consider the strongly connected subtournament $T(W)$
induced by $W$.
Since $T(W)$ contains $T_{ncr},$ but $T_{ncr}$ is not contained
in any proper strongly connected subtournament of $T,$
we have $T(W)=T$.
In particular, this means that $Q$ includes all critical vertices of $T.$
The condition
$|T_{cr}|\ge \ell-1$
implies that $|Q|\ge |T_{cr}|+2\ge \ell+1.$
Hence, $c_{\ell}\bigl(T(Q)\bigr)=c_{\ell}\bigl(T_{|Q|-1,|Q|}\bigr)=
|Q|-\ell+1.$
Moreover,  the set $\overline{Q}$
of vertices not belonging to $Q$
is a subset of $V(T_{ncr})$ and hence,
by our assumption, $c_{\ell}(T,w)=1$ for any vertex $w\in \overline{Q}.$
This implies that
$$c_{\ell}(T)\le c_{\ell}\bigl(T(Q)\bigr)+\sum\limits_{w\in \overline{Q}}
c_{\ell}(T,w)=|Q|-\ell+1+(n-|Q|)=n-\ell+1.$$
By Theorems II and III,
we have $T\cong T_{n-1,n}.$
The theorem is proved. $\blacksquare$

\smallskip

Let $\Cal{T}_{\le d,n}$ be the class
of all strongly connected tournaments
of order $n$ whose diameter is not greater than $d.$
Theorem 1 implies that for any $T\in {\Cal T}_{\le n-2,n},$
there exists a non-critical vertex $w$ with $c_{\ell}(T,w)\ge 2.$
What is the minimum possible value of $c_{\ell}(T,w),$
where $w\in V(T_{ncr}),$
that is necessary for proving Conjecture 1 in
the case of arbitrary parameter $d\ge \frac{n+3}{2}$?

To answer this natural question, let us consider the non-critical vertices
of $T_{d,n}^{-}.$ It is not difficult to check that
for this case, $V(T_{cr})=\{v_{1},...,v_{d-1}\}$
(see Lemma 2).
Choose any vertex $w$ of the subtournament $TT_{\lceil\frac{n-d+1}{2}\rceil}$
in $T_{d,n}^{-}$ and then
take any subset $S$ of order $0\le k\le \lceil\frac{n-d+1}{2}\rceil-1$ in $TT_{\lceil\frac{n-d+1}{2}\rceil}-w.$
Let $w_{1},...,w_{\lceil\frac{n-d+1}{2}\rceil}$
be the unique hamiltonian path
of $TT_{\lceil\frac{n-d+1}{2}\rceil}.$
Then $S=\{w_{i_{1}},...,w_{i_{k}}\},$
where $i_{1}<...<i_{k}.$
Obviously, $w=w_{p},$ where
$i_{s}<p<i_{s+1}$
for some $0\le s\le k.$
(Here $i_{0}=0$ and $i_{k+1}=\lceil\frac{n-d+1}{2}\rceil+1$.)
Moreover, $k+3\le \lceil\frac{n-d+1}{2}\rceil+2
\le n-d+3.$
Hence, if
$n-d+3\le \ell\le d,$ then
$w_{i_{1}},...,w_{i_{s}},w_{p},w_{i_{s+1}},...,w_{i_{k}},
v_{1},...,v_{\ell-k-1},w_{i_{1}}$ is a circuit of length $\ell.$
It is easy to see that any circuit of length $\ell$
including $w$ admits this form. In other words,
under the restrictions
on $\ell$ given in the assumptions of Conjecture 1,
there is a one-to-one
correspondence between the circuits of length $\ell$ passing through $w$
and the subsets of the vertex-set of the tournament
$TT_{\lceil\frac{n-d+1}{2}\rceil}-w$ of order $\lceil\frac{n-d+1}{2}\rceil-1.$
Hence, we have
$$c_{\ell}(T_{d,n}^{-},w)=\sum\limits_{k=0}^{\lceil\frac{n-d+1}{2}\rceil-1}
\binom{\lceil\frac{n-d+1}{2}\rceil-1}{k}=
2^{\lceil\frac{n-d+1}{2}\rceil-1}.$$
It is not difficult to check that if $d\le n-2,$ then
$T_{d,n}^{-}-w\cong T_{d,n-1}.$
So, to prove Conjecture 1 by induction on $n,$ it would not be bad to show that
the following suggestion holds.

\smallskip

{\bf Conjecture 2.} {\sl Let $T$ be a strongly connected tournament of order $n\ge 4$
whose diameter does not exceed $d,$ where $d\ge \frac{n+3}{2}.$
Then for each $\ell=n-d+3,...,d,$ there exists a non-critical vertex $w$ of $T$
such that the diameter of the strongly connected subtournament $T-w$ is not
greater than $d,$ either, and}
$$c_{\ell}(T,w)\ge 2^{\lceil\frac{n-d+1}{2}\rceil-1}.$$

\smallskip

In the sequel, we also need a sharp lower bound on
$|T_{ncr}|$ in the class ${\Cal T}_{d,n}.$
We start its presentation with the case $d=2.$

\smallskip

{\bf Proposition 1 [16].} {\sl A strongly connected
tournament $T$ with diameter two
contains at most $3$ critical vertices. Let $1\le c\le 3.$ Then
$T$ includes exactly $c$ critical vertices if and only if
$T$ is obtained from the cyclic triple $\Delta$ by replacing
precisely $3-c$ of its vertices by
strongly connected tournaments of order at least $3$ with diameter $2.$}

\smallskip

{\bf Remark 2}. Proposition 1
taken together with Bush's revisions [4] of
Moon's theorem on hamiltonian paths in strongly connected tournaments
[10] and
Thomassen's theorem on hamiltonian circuits in $2$-strong tournaments [18]
yields a lower bound on the number of circuits of length $n-h$
in the class ${\Cal T}_{2,n}$ (see Theorem 4 [16]).
In particular, this bound implies that for a given $h\ge 0,$
the minimum of $c_{n-h}(T)$ in the class ${\Cal T}_{2,n}$
grows exponentially with respect to $n.$
Note that for given $d\ge 3$ and $h\ge 0,$
this minimum in the class ${\Cal T}_{d,n}$ is $O(n^{h}).$

Obviously, Proposition 1 implies that any strongly connected tournament
of order $n\ge 4$ with diameter $2$ contains at least $n-2$ non-critical vertices.
For $d\ge 3,$ a sharp lower bound on $|T_{ncr}|$
including $n$ and $d$ is also known
and was first obtained in [16] with the use of Lemma 1.

\smallskip

{\bf Proposition 2 [16].}
{\sl Let $T$ be a strongly connected tournament
of order $n\ge 4$ with diameter $d\ge 3.$
Then $T$
contains at least $n-d+1$ non-critical vertices.
Moreover, if $|T_{ncr}|=n-d+1,$ then $T_{cr}\cong T_{d-2,d-1}$.}

\smallskip

{\bf Remark 3}.
Proposition 2 and Camion's theorem (see Theorem II)
imply that $c_{n-1}(T)\ge n-d+1$
for any $T\in \Cal{T}_{d,n},$ where $n>d\ge 3.$
Moreover, based on Proposition 2, one can describe
all strongly connected tournaments in the class $\Cal{T}_{d,n}$
admitting exactly
$n-d+1$ circuits of length $n-1.$ This has been recently done
in [16] (see Theorem 2).  In particular, this theorem implies that
any tournament $T\in \Cal{T}_{d,n}$ with $c_{n-1}(T)=n-d+1,$
where $n\ge 6,$
has exactly one hamiltonian circuit and hence, can be simply described
in terms of Douglas' theorem [6] (see also Proposition 3 below).

\smallskip

Unfortunately, at the moment,
we are unable to classify (in a simple form)
all strongly connected tournaments in $\Cal{T}_{d,n}$ containing exactly
$n-d+1$ non-critical vertices. However, we can do this perfectly if
$|T_{ncr}|=n-d+1$ and
$T_{ncr}$ contains precisely two strongly connected components.
(Note that if $T_{ncr}$ is strong and $T\neq T_{ncr},$ then
by Lemma 1, $T=\Delta(v_{1},v_{2},T_{ncr});$ in particular, for this
case, $d\le 3$; so, if $d\ge 4$ and $T$ is not $2$-strong,
then $T_{ncr}$ has at least two strong components.)
In some sense, the following proposition can be considered
as a natural analog of Theorem V for the class ${\Cal T}_{d,n}.$

\smallskip

{\bf Corollary 1.}
{\sl Let $T$ be a strongly connected tournament of order $n\ge 4$
with diameter $d\ge 3.$ Assume that $T$
admits exactly $n-d+1$ non-critical vertices (so, by Proposition 2,
we can assume that $T_{cr}=T_{d-2,d-1}(v_{1},...,v_{d-1})$)
and $T_{ncr}$ contains precisely two strongly connected components
$T_{1}$ and $T_{2}$, i.e. $T_{ncr}=TT_{2}(T_{1},T_{2}).$ Then
$T=T_{d,d+1}(T_{2},v_{1},...,v_{d-1},T_{1})$.}

\smallskip

{\sl Proof.}  Assume first that $|T_{ncr}|=n-d+1$
and the number $m$ of strongly connected components $T_{1},...,T_{m}$
of $T_{ncr}$ is not less than $2.$
Without loss of generality, we can assume that
$T_{ncr}=TT_{m}(T_{1},...,T_{m}).$
Consider a shortest path $Q=v_{0},v_{1},...,v_{p-1},v_{p}$ from $T_{m}$
to $T_{1}$ in $T$.
Obviously, $T(Q)=T_{p,p+1}(v_{0},...,v_{p})$
and $p\le d.$ Lemma 1 implies that $Q$ includes all critical vertices of
the strong tournament $T.$
By the assumption, $|T_{cr}|=d-1$ and hence,
in fact, $p=d.$ In particular, this means that
$T_{cr}=T(Q)-v_{0}-v_{d}=T_{d-2,d-1}(v_{1},...,v_{d-1}).$

Let us now determine the orientation of the arcs
between $T_{m}$ and $T_{cr}.$
Obviously, if $w\to v_{i},$ where $w\in V(T_{m})$ and $2\le i\le d-1,$
then the path $w,v_{i},...,v_{d-1},v_{d}$ from $T_{m}$ to $T_{1}$
is shorter than $Q.$ Hence, $\{v_{2},...,v_{d-1}\}
\Rightarrow T_{m}.$ Suppose now that
$(w,v_{1})$ is not an arc for some $w\in V(T_{m}).$
Since $\{T_{1},...,T_{m-1},v_{2},...,v_{d-1}\}
\Rightarrow T_{m},$ a shortest path from $w$ to $v_{d}$ must have the form
$w,w_{1},...,w_{k},v_{1},...,v_{d},$ where $k\ge 1$
and $w_{1},...,w_{k}$ is a path in $T_{m},$ and hence,
the distance between $w$ and $v_{d}$ is strictly greater than $d,$
which is impossible.
This means that $T_{m}\Rightarrow v_{1}.$ Similarly,
$v_{d-1}\Rightarrow T_{1}\Rightarrow\{v_{1},...,v_{d-2}\}.$
In particular, if $m=2,$ then $T=T_{d,d+1}(T_{m},v_{1},...,v_{d-1},T_{1}).$
The corollary is proved. $\blacksquare$

\smallskip

Lemma 2 implies that in fact,
any strongly connected tournament of the form $T_{d,d+1}(T_{2},v_{1},...,v_{d-1},T_{1})$ with
$|T_{1}|+|T_{2}|=n-d+1$ contains exactly $n-d+1$ non-critical vertices.
Note that for achieving this (sharp) lower bound on $|T_{ncr}|$
in the class ${\Cal T}_{d,n},$ where $n>d\ge 3,$
the tournaments $T_{1}$ and $T_{2}$
in $T_{d,d+1}(T_{2},v_{1},...,v_{d-1},T_{1})$
need not be strongly connected. According to our Conjecture 1, for
$d\ge \frac{n+3}{2}$ and
$\ell=n-d+3,...,d,$ the lower bound on $c_{\ell}(T)$
in the class $\Cal{T}_{\le d,n}$
is attained if and only if
either of $T_{1}$ and $T_{2}$ is transitive and their orders
are as nearly equal as possible.
In the next section, we reprove this statement for the case $d=n-2$
and give an original proof for $d=n-3.$

\bigskip

\centerline{\bf 3. The proof of the statement of Conjecture $\bold{1}$ for $\bold{d=n-2}$ and $\bold{d=n-3}$}

\smallskip

To prove the statement of Conjecture $1$
for $d=n-3$ by induction, we also need the fact that
it holds for $d=n-2.$
As we have noted above, the case $d=n-2$ was completely considered
in our paper [16].
Unfortunately, the proof given in [16]
is complicated enough.  For the convenience of the reader, below
we give an easier proof which is based on Theorem 1.

Let
$T_{n-2,n}(z_{0},...,\widehat{z_{n-3},z_{n-2},z_{n-1}})$ and
$T_{n-2,n}^{-}(\widehat{z_{0},z_{1},z_{2}},...,z_{n-1})$ be the
result of the reversion of the arcs $(z_{n-1},z_{n-3})$ and
$(z_{2},z_{0}),$ respectively, in $T_{n-1,n}(z_{0},z_{1},z_{2},$
$...,z_{n-3},z_{n-2},z_{n-1}).$
So, $\widehat{z_{0},z_{1},z_{2}}$ and $\widehat{z_{n-3},z_{n-2},z_{n-1}}$
mean that the subtournaments of order $3$ induced by the vertex-sets
$\{z_{0},z_{1},z_{2}\}$ and $\{z_{n-3},z_{n-2},z_{n-1}\}$ are transitive.
Obviously, $T_{n-2,n}$ and $T_{n-2,n}^{-}$ given above
are obtained from $T_{n-2,n-1}(z_{0}$
$,...,z_{n-2})$
and $T_{n-2,n-1}(z_{1},...,z_{n-1})$ by
replacing the right and left non-critical vertices
$z_{n-2}$ and $z_{1}$ with the transitive tournaments
$TT_{2}(z_{n-2},z_{n-1})$
and $TT_{2}(z_{0},z_{1}),$ respectively, and hence, are isomorphic to
the tournaments $T_{n-2,n}$ and $T_{n-2,n}^{-}$ introduced
before the statement of Conjecture 1. One can say that
$z_{n-1}$ and $z_{1}$ enlarge $z_{n-2}$ and $z_{0},$ respectively.
In the sequel, we call $z_{0}$ and $z_{n-1}$  the left and right
non-critical vertices in
$T_{n-2,n}(z_{0},...,\widehat{z_{n-3},z_{n-2},z_{n-1}})$ and
$T_{n-2,n}^{-}(\widehat{z_{0},z_{1},z_{2}},...,z_{n-1}),$
respectively.

Recall that $T_{n-1,n}$ is the unique strong tournament of order $n$
with a pair of vertices at distance $n-1$ from each other.
Strongly connected tournaments $T_{n-2,n}$
and $T_{n-2,n}^{-}$ can be similarly defined.
More precisely, $T_{n-2,n}$ is uniquely characterized by the property
that there exists a vertex at distance $n-2$ to two other vertices
and $T_{n-2,n}^{-}$ is uniquely determined by
the condition that there exists a vertex at distance $n-2$ from two other
vertices. In both cases, all the three vertices
(and only they) are non-critical.

\smallskip

{\bf Theorem 2 [16].}
{\sl Let
$T$ be a strongly connected tournament
of order $n\ge 7$ whose diameter is not greater than
$n-2.$ (In other words, $T$ is not isomorphic to $T_{n-1,n}$.)
Then $c_{\ell}(T)\ge n-\ell+2$
for any given $\ell$ in the range $5,...,n-2$ with equality holding iff
$T$ is isomorphic to $T_{n-2,n}$ or $T_{n-2,n}^{-}$.}

{\sl Proof.}
By Theorem III,
if $T\in \Cal{T}_{\le n-2,n},$ then $c_{\ell}(T)\ge n-\ell+2$
for any $\ell=4,...,n-1.$
Assume that
$c_{\ell}(T)=n-\ell+2$
for some $\ell=5,...,n-2.$
By Theorem 1,
we can choose a non-critical vertex $w$ with $c_{\ell}(T,w)\ge 2.$
For this case, $c_{\ell}(T-w)\le n-\ell.$
On the other hand, by Theorem II, we have $c_{\ell}(T-w)\ge (n-1)-\ell+1=n-\ell.$
This means that $c_{\ell}(T,w)=2$ and
$c_{\ell}(T-w)=(n-1)-\ell+1.$ The last equality taken together with
Theorem III implies that $T-w\cong T_{n-2,n-1}.$
More precisely, in the sequel, we assume that
$T-w=T_{n-2,n-1}(z_{0},...,z_{n-2}).$

Since $T$ is strongly connected,
there always exists $k$ satisfying the inequality $0\le k\le n-3$
such that either $z_{k}\to w\to z_{k+1}$ or
$\{z_{k+1},...,z_{n-2}\}\Rightarrow w\Rightarrow\{z_{0},...,z_{k}\}.$
Below we consider these two possibilities separately.
Note that the cases  $k=n-4$ and $k=n-3$ are reduced to the cases
$k=1$ and $k=0,$ respectively,
by the transition from the original tournament $T$ to its converse.
Hence, in the sequel,
we can assume that $0\le k\le n-5.$

{\bf (I)} Suppose first
that $z_{k}\to w\to z_{k+1}$ for some $0\le k\le n-5.$

{\bf (1)} Case $0\le k\le \ell-3.$ In the next three subitems,
we show that it is possible only if
$T-z_{0}=T_{n-2,n-1}(w,z_{1},...,z_{n-2}).$

{\bf (1a)} Let us prove first that the subtournament
$T-z_{0}$ is strongly connected. Obviously, it is so for the
subtournament $T-w-z_{0}$ which also coincides with $T-z_{0}-w.$
It is known (and easy to check) that adding to a strongly connected tournament
a new vertex $w$ which dominates at least one of its vertices and is dominated
by at least one of its vertices leads to a strongly connected digraph, again.
Hence, if $T-z_{0}$ is not strongly connected,
then either $w\Rightarrow\{z_{1},...,z_{n-2}\}$
or $\{z_{1},...,z_{n-2}\}\Rightarrow w.$ By condition
$z_{k}\to w\to z_{k+1}$,
this is possible if and only if $z_{0}\to w\Rightarrow\{z_{1},...,z_{n-2}\}.$
However, this contradicts
the equality $c_{\ell}(T,w)=2$ established above because

$(*)$ if $z_{0}\to w\Rightarrow\{z_{1},...,z_{n-3}\},$ then
$\gamma_{1}=z_{0},w,z_{1},...,z_{\ell-2},z_{0},$\
$\gamma_{2}=z_{0},w,z_{2},...,$ $z_{\ell-1},z_{0}$ and
$\gamma_{3}=z_{0},w,z_{3},...,z_{\ell},z_{0}$
will be three circuits of length $\ell$
including the vertex $w.$

{\bf (1b)} Note that $z_{0},...,z_{k},w,z_{k+1},...,z_{\ell-2},z_{0}$
and $z_{0},z_{1},...,z_{\ell-1},z_{0}$
are two circuits of length $\ell$
passing through the vertex $z_{0}.$
Thus, $c_{\ell}(T-z_{0})=(n-1)-\ell+1$ and hence,
by Theorem III, we have $T-z_{0}\cong T_{n-2,n-1}.$

{\bf (1c)} Obviously,
$T-z_{0}-w=T_{n-3,n-2}(z_{1},...,z_{n-2}).$
If $w$ is the right non-critical vertex of $T-z_{0},$
then by the structure  of $T_{n-2,n-1},$
it is dominated by the right non-critical vertex of $T-z_{0}-w$
(i.e. $z_{n-2}$) and dominates all the other vertices of $T-z_{0}-w.$
Since $z_{k}\to w\to z_{k+1}$ for some $0\le k\le \ell-3\le n-5,$
we have
$z_{0}\to w \Rightarrow \{z_{1},...,z_{n-3}\},$ which
is impossible by virtue of $(*)$.
Therefore, $w$ is the left non-critical
vertex of $T-z_{0}$.
By the structure  of $T_{n-2,n-1},$
it dominates the left non-critical vertex of $T-z_{0}-w$ (i.e. $z_{1}$)
and is dominated
by all the other vertices of $T-z_{0}-w.$
Thus,
$T-z_{0}=T_{n-2,n-1}(w,z_{1},...,z_{n-2}).$

{\bf It remains} to note that if
$z_{0}\to w,$ then
$T=T_{n-2,n}^{-}(\widehat{z_{0},w,z_{1}},...,z_{n-2}).$
Otherwise (i.e. in the case of $w\to z_{0}$), we have
$\{z_{2},...,z_{n-2}\}\Rightarrow w\Rightarrow\{z_{0},z_{1}\},$
which is inconsistent with the existence of $k$
for which $z_{k}\to w\to z_{k+1}.$

{\bf (2)}  Case $\ell-2\le k\le n-5.$
In this situation, the vertex $w$ is contained in three circuits
$\gamma_{1}=z_{k-\ell+3},...,z_{k},w,z_{k+1},z_{k-\ell+3},$\
$\gamma_{2}=z_{k-\ell+4},...,z_{k},w,z_{k+1},z_{k+2},z_{k-\ell+4},$ and
$\gamma_{3}=z_{k-\ell+5},...,z_{k},w,z_{k+1},z_{k+2},z_{k+3},z_{k-\ell+5}$
of length $\ell,$ which is impossible.

{\bf (II)}
Suppose now that $\{z_{k+1},...,z_{n-2}\}\Rightarrow w
\Rightarrow\{z_{0},...,z_{k}\}$ for some $0\le k\le n-5.$
(Obviously, this condition taken together with
$T-w=T_{n-2,n-1}(z_{0},...,z_{n-2})$ uniquely determines
$T.$) 

{\bf (1)} Case $0\le k\le 1.$
If $k=0,$ then $T=T_{n-1,n}(w,z_{0},...,z_{n-2})$
and hence, the diameter of $T$ equals $n-1,$ which is impossible.
In turn, if $k=1,$ then $T=T_{n-2,n}^{-}(\widehat{w,z_{0},z_{1}},...,z_{n-2})$.

{\bf (2)} Case $2\le k\le \ell-3.$ Then
the vertex $w$ is contained in three circuits
$\gamma_{1}=w,z_{0},...,$ $z_{\ell-2},w,$
$\gamma_{2}=w,z_{1},...,z_{\ell-1},w,$ and
$\gamma_{3}=w,z_{2},...,z_{\ell},w$ of length $\ell,$
which is impossible.

{\bf (3)} Case $\ell-2\le k \le n-5.$
Then the vertex $w$ belongs to three circuits
$\gamma_{1}=w,z_{k-\ell+3},...,z_{k+1},w,$
$\gamma_{2}=w,z_{k-\ell+4},...,z_{k+2},w,$ and
$\gamma_{3}=w,z_{k-\ell+5},...,z_{k+3},w$ of length $\ell,$ which is impossible.
The theorem is proved. $\blacksquare$

\smallskip

A simple modification of
the proof of Theorem 2 also allows us to consider the case $\ell=4.$
Indeed, it is easy to see that we need to explain only the case
$2\le k\le n-5.$
If $z_{k}\to w\to z_{k+1},$ then
for $\ell=4,$
the cycle
$\gamma_{3}=z_{k-\ell+5},...,z_{k},w,z_{k+1},z_{k+2},...,$ $z_{k+3},z_{k-\ell+5}$
makes no sense.
To present a new $\gamma_{3},$ consider
the subtournament $S_{4}$ induced by the vertex-set
$V_{4}=\{z_{k-2},z_{k-1},z_{k},w\}.$
Obviously, the vertex-set
$\{z_{k-2},z_{k-1},z_{k}\}$ induces a cyclic triple.
Hence, if $S_{4}$
is not strongly connected, then $\{z_{k-2},z_{k-1},z_{k}\}\Rightarrow w$
and hence, $\gamma_{3}=z_{k-2},w,z_{k+1},z_{k+2},z_{k-2}$
is a third circuit of length $4$ containing the vertex $w.$
In the opposite case, we can take
a hamiltonian circuit of
the strongly connected subtournament $S_{4}$ of order $4$ as $\gamma_{3}.$
In turn, if
$\{z_{k+1},...,z_{n-2}\}\Rightarrow w\Rightarrow\{z_{0},...,z_{k}\},$
then for $\ell=4,$ the cycle $\gamma_{3}=w,z_{k-\ell+5},...,z_{k+3},w$
is impossible.
For this case, there exist exactly two cycles of length $\ell=4$
containing the vertex $w,$ namely,
$\gamma_{1}=w,z_{k-\ell+3},...,z_{k+1},w$ and
$\gamma_{2}=w,z_{k-\ell+4},...,z_{k+2},w$
(note that the properties of the subtournament
$T-w=T_{n-2,n-1}(z_{0},...,z_{n-2})$ imply that
any path of length $2$
from $N^{+}(w)=\{z_{0},...,z_{k}\}$ to $N^{-}(w)=\{z_{k+1},...,z_{n-2}\}$
is either $z_{k-1},z_{k},z_{k+1}$ or $z_{k},z_{k+1},z_{k+2}$).
For such $T,$ we get $c_{4}(T)=c_{4}(T-w)+c_{4}(T,w)=
(n-1)-4+1+2=n-2.$
Hence, the equality $c_{4}(T)=n-2$ holds if and only if
$T$ is isomorphic to the tournament obtained from
$T_{n-2,n-1}(z_{0},...,z_{n-2})$ by adding a vertex $w$ with
$\{z_{k+1},...,z_{n-2}\}\Rightarrow w\Rightarrow\{z_{0},...,z_{k}\},$
where $1\le k \le n-4.$
This class contains $n-4$ elements and coincides with
the class ${\Cal H}_{n-2,n}$ of all strongly connected tournaments
of order $n$
and diameter $n-2$ admitting exactly one hamiltonian circuit and
including precisely three non-critical vertices.
For arbitrary $d\ge 3,$ the class ${\Cal H}_{d,n}$
will be introduced and described in the next section (see Proposition 3).
By [3], any $T\in {\Cal H}_{n-2,n}$ also has the least possible
number of cycles of length $3.$
Note that the equality $c_{n-1}(T-w)=1$ does not mean that
$T-w\cong T_{n-2,n-1}.$ Hence, the proof of
Theorem 2 cannot be adapted to the case $\ell=n-1.$
Nevertheless, according to Theorem 2 [16], for $n\ge 7,$
the equality $c_{n-1}(T)=3$ also holds if and only if
$T\in {\Cal H}_{n-2,n}.$

\smallskip

Let
$T_{n-3,n}(\widehat{z_{0},z_{1},z_{2}},...,\widehat{z_{n-3},z_{n-2},z_{n-1}})$ be the result of the
reversion of the arcs $(z_{2},z_{0})$ and $(z_{n-1},z_{n-3})$ in
$T_{n-1,n}(z_{0},...,z_{n-1}).$
Obviously, it is obtained from $T_{n-3,n-2}(z_{1},...,z_{n-2})$
by replacing its right and left non-critical vertices
$z_{n-2}$ and $z_{1}$ with the transitive tournaments
$TT_{2}(z_{n-2},z_{n-1})$
and $TT_{2}(z_{0},z_{1}),$ respectively, and hence, is isomorphic to
the tournament $T_{n-3,n}$ introduced in Section 1.

\smallskip

{\bf Theorem 3.}
{\sl Let
$T$ be a strongly connected tournament of order $n\ge 9$
whose diameter is not greater than
$n-3.$
Then $c_{\ell}(T)\ge n-\ell+3$
for any given $\ell$ in the range $6,...,n-3$ with equality holding iff
$T$ is isomorphic to $T_{n-3,n}$.}

\smallskip

{\sl Proof.}
By Theorem III,
if $T\in \Cal{T}_{\le n-1,n}$  and
$c_{\ell}(T)=n-\ell+1$
for some $\ell=4,...,n-1,$ then $T\cong T_{n-1,n}.$
In turn, by Theorem 2,
if $c_{\ell}(T)=n-\ell+2$
for some $\ell=5,...,n-2,$ then $T\in \Cal{T}_{n-2,n}.$
Hence, if $T\in \Cal{T}_{\le n-3,n},$ then $c_{\ell}(T)\ge n-\ell+3$
for any $\ell=5,...,n-2.$

Let us now show that if $T\in \Cal{T}_{\le n-3,n}$ and
$c_{\ell}(T)=n-\ell+3$
for some $\ell=6,...,n-3,$ then
$T\cong T_{n-3,n}$.
Let $w$ be a non-critical vertex with $c_{\ell}(T,w)\ge 2$
(recall that we can always take such a vertex by virtue of Theorem 1).
Since $c_{\ell}(T)=n-\ell+3,$ we have $c_{\ell}(T-w)\le n-\ell+1.$
On the other hand,
by Theorem II, $c_{\ell}(T-w)\ge n-\ell$ and hence,
$c_{\ell}(T,w)\le 3.$
Rewrite the bounds obtained for $c_{\ell}(T-w)$ as
$(n-1)-\ell+1\le c_{\ell}(T-w)\le (n-1)-\ell+2.$
By Theorem III and Theorem 2,
we can assume that the subtournament
$T-w$ of order $n-1$
is obtained from
$T_{n-2,n-1}(z_{0},...,z_{n-2})$
by reversing at most one of the arcs $(z_{n-2},z_{n-4})$
and $(z_{2},z_{0}).$

Since the cases $k=n-5,n-4,n-3$ are reduced to
the cases $k=2,1,0,$ respectively,
by means of the transition
from the original tournament $T$ to its converse,
without loss of generality,
we can assume that $0\le k\le n-6.$
Suppose first that there exists $k$ such
that $z_{k}\to w\to z_{k+1}.$
If $\ell-2\le k\le n-6,$ then the vertex $w$ is contained in four circuits
$\gamma_{1},$ $\gamma_{2},$
$\gamma_{3},$ and $\gamma_{4},$ where
$\gamma_{i}=z_{k-\ell+2+i},...,z_{k},w,z_{k+1},...,z_{k+i},z_{k-\ell+2+i},$
of length $\ell,$ which is impossible.
So, it suffices to consider only the case $0\le k\le \ell-3.$
Repeating the arguments presented in subitem {\bf (1a)}
of the proof of Theorem 2
we get that the subtournament $T-z_{0}$
is also strongly connected. (It is not difficult to check
that such argumentation
does not depend on the orientation of the arcs between the vertices
$z_{0}$ and $z_{2}$ and between the vertices $z_{n-4}$ and $z_{n-2}$).
Only now one
should replace the statement $(*)$ therein by
the fact that

(**) if $z_{0}\to w\Rightarrow\{z_{1},...,z_{n-5}\},$
then
$\gamma_{1}=z_{0},w,z_{1},...,$ $z_{\ell-2},z_{0},$\
$\gamma_{2}=z_{0},w,z_{2},...,$
$z_{\ell-1},z_{0},$\
$\gamma_{3}=z_{0},w,z_{3},...,z_{\ell},z_{0},$
and $\gamma_{4}=z_{0},w,z_{4},...,z_{\ell+1},z_{0}$
are four circuits of length $\ell$
including the vertex $w,$
which is impossible by virtue of $c_{\ell}(T,w)\le 3.$

In turn, subitem {\bf (1b)} of the proof of Theorem 2
implies that $c_{\ell}(T,z_{0})\ge 2.$ Thus,
$(n-1)-\ell+1\le c_{\ell}(T-z_{0})\le (n-1)-\ell+2.$
By Theorem III, the equality $c_{\ell}(T-z_{0})=(n-1)-\ell+1$
is possible only if
$T-z_{0}\cong T_{n-2,n-1}.$
More precisely, subitem {\bf (1c)}
of the proof of Theorem 2 allows us to state that
$T-z_{0}=T_{n-2,n-1}(w,z_{1},...,z_{n-2}).$
However, for this case, regardless of the orientation of the arc
between the vertices $z_{0}$ and $z_{2},$
the distance between $w$ and $z_{n-2}$ in $T$ itself
is also equal to $n-2,$ which is impossible.
Hence,
$c_{\ell}(T-z_{0})=(n-1)-\ell+2.$
The last equality and
Theorem 2 mean that either $T-z_{0}\cong T_{n-3,n-1},$
or $T-z_{0}\cong T_{n-3,n-1}^{-}$
(their structure is described in detail in the paragraph placed
just before the statement of Theorem 2).
Note also that if
$z_{0}\to z_{2},$ then the sequence
$z_{0},z_{2},...,z_{\ell},z_{0}$ forms a third circuit
of length $\ell$ containing $z_{0}$ and hence,
$c_{\ell}(T-z_{0})=(n-1)-\ell+1,$ which, as we have seen above,
is impossible. Thus, $z_{2}\to z_{0}$ and so, only the orientation
of the arc betweeen $z_{n-4}$ and $z_{n-2}$ can be arbitrary.

Assume first that $z_{n-4}\to z_{n-2}.$ Then we have
$T-z_{0}-w=T_{n-4,n-2}
(z_{1},...,$
$\widehat{z_{n-4},z_{n-3},z_{n-2}}).$
Since $T-z_{0}$ is isomorphic to $T_{n-3,n-1}$ or $T_{n-3,n-1}^{-}$
and $T-z_{0}-w\cong T_{n-4,n-2},$
we have $T-z_{0}\cong T_{n-3,n-1}$ and $w$ is the
left non-critical vertex of $T-z_{0}.$ By the structure of $T_{n-3,n-1},$
it dominates the left non-critical vertex of $T-z_{0}-w$
(in our notation, this is $z_{1}$) and is dominated by
all the other vertices of $T-z_{0}-w.$
Thus,
$T-z_{0}=T_{n-3,n-1}(w,z_{1},...,$
$\widehat{z_{n-4},z_{n-3},z_{n-2}})$.
If $z_{0}\to w,$ then
$T=T_{n-3,n}(\widehat{z_{0},w,z_{1}},...,\widehat{z_{n-4},z_{n-3},z_{n-2}}),$
which is completely in agreement with the statement of the theorem.
Otherwise, we have
$\{z_{2},...,z_{n-2}\}\Rightarrow w\Rightarrow\{z_{0},z_{1}\},$
which is inconsistent with the existence of $k$ with $z_{k}\to w\to z_{k+1}.$

Assume now that $z_{n-2}\to z_{n-4}$ and hence,
$T-z_{0}-w=T_{n-3,n-2}(z_{1},...,z_{n-2}).$
If either $T_{n-3,n-1}-w\cong T_{n-3,n-2}$ or $T_{n-3,n-1}^{-}-w\cong
T_{n-3,n-2},$ then  the vertex
$w$ enlarges either the right non-critical vertex of $T_{n-3,n-2}$
(in our notation, this is $z_{n-2}$) or the left non-critical vertex of
$T_{n-3,n-2}$ (in our notation, this is $z_{1}$).
For the first case,
$w\Rightarrow\{z_{1},...,z_{n-4}\}$ (this can hold only for $k=0$),
which is impossible by virtue of $(**).$ For the second case, we have
$\{z_{3},...,z_{n-2}\}\Rightarrow w\to z_{2}$ and hence, the distance
between $z_{0}$ and $z_{n-2}$ in $T$ is equal to $n-2$ (recall that
$z_{2}\to z_{0}$ if $z_{k}\to w\to z_{k+1},$ where
$0\le k\le \ell-3).$
However, this is in contradiction with the assumption of the theorem.

{\bf Finally}, suppose that there exists $k$ such that
$\{z_{k+1},...,z_{n-2}\}\Rightarrow w\Rightarrow\{z_{0},...,$
$z_{k}\}.$
Recall that the subtournament $T-w$ is obtained
from $T_{n-2,n-1}(z_{0},...,z_{n-2})$ by reversing
at most one of the arcs $(z_{n-2},z_{n-4})$
and $(z_{2},z_{0}).$
Obviously, if $T-w\cong T_{n-2,n-1},$
then the distance between $z_{0}$ and $z_{n-2}$ in $T$
equals $n-2.$ Hence, either
$T-w=T_{n-3,n-1}(z_{0},...,\widehat{z_{n-4},z_{n-3},z_{n-2}})$ or
$T-w=T_{n-3,n-1}^{-}(\widehat{z_{0},z_{1},z_{2}},...,z_{n-2})$
and, in both cases, $c_{\ell}(T,w)=2.$
If $0\le k\le 1,$ then either
the distance between $w$ and $z_{n-2}$ in $T$
is equal to $n-2,$ which is impossible, or we have
$T=$
$T_{n-3,n}(\widehat{w,z_{0},z_{1}},...,\widehat{z_{n-4},z_{n-3},z_{n-2}})$.
The latter holds iff $k=1$ and
$T-w=T_{n-3,n-1}(z_{0},...,\widehat{z_{n-4},z_{n-3},z_{n-2}}).$
For the cases $2\le k\le \ell-3$ and $\ell-2\le k\le n-6,$
a contradiction with the equality $c_{\ell}(T,w)=2$
is established absolutely in the same way as it has been done
in items {\bf (2)} and {\bf (3)} of section II of the proof
of Theorem 2, respectively. The theorem is proved. $\blacksquare$

\bigskip

\centerline{\bf 4. Strongly connected tournaments
of order $\bold{n}$ and diameter $\bold{\le n-3}$}

\smallskip

\centerline{\bf with exactly $\bold{c_{\pmb\ell}(T_{n-3,n})}$
circuits of length $\bold{\pmb\ell}$ for the remaining values of $\pmb\ell$}

\smallskip

Is the restriction on the length
$\ell$ in the condition of Theorem 3 essential?
To answer this natural question, let us consider
the class ${\Cal H}_{d,n}$ of all strongly connected
tournaments of order $n$ with diameter $d$
that admit exactly one hamiltonian
circuit and contain exactly $n-d+1$ non-critical vertices,
where $n>d\ge 3.$
The following proposition describes the representatives of
the class ${\Cal H}_{d,n}$
in terms of parameters of Douglas' theorem [6].

\smallskip

{\bf Proposition 3 [16]}.
{\sl A strongly connected tournament $T$ of order $n\ge 4$

with diameter $d\ge 3$ is contained in the class ${\Cal H}_{d,n}$ if and only if

the following two conditions hold:

(1) there exist two
vertex-sets $V=\{v_{1},...,v_{d-1}\}$ and $W=\{w_{1},...,w_{n-d+1}\}$ such

that

(a) $T(V)=T_{d-2,d-1}(v_{1},...,v_{d-1}),$

(b) $T(W)=TT_{n-d+1}(w_{1},...,w_{n-d+1});$

(2) for each $i=1,...,n-d+1,$ there exists $h_{i}$
such that
$\{v_{h_{i}+1},...,v_{d-1}\}\Rightarrow w_{i}\Rightarrow$

$\{v_{1},...,v_{h_{i}}\}$ with

(a) $h_{1}=d-2,$ $1\le h_{i}\le d-2$ for $i=2,...,n-d,$ and $h_{n-d+1}=1,$

(b) $h_{s}\le h_{r}+1\ \text{ for any }\ 1\le r\le s \le n-d+1.$}

\smallskip

A walk of a tournament is called {\sl spanning} if
it contains all the vertices. Note that a tournament is
strongly connected if and only if it admits at least one
spanning closed walk. Let us discuss the conditions of Proposition 3.
Obviously, the closed walk $w_{i+1},...,w_{n-d+1},v_{1},...,v_{d-1},
w_{1},...,w_{i-1},w_{i+1}$ is spanning in $T-w_{i}.$ Hence, for each
$i=1,...,n-d+1,$ this subtournament is strongly connected. No arc leaves
the vertex $w_{n-d+1}$ in $T-v_{1}$ and no arc enters the vertex $w_{1}$
in $T-v_{d-1}.$ Hence, they are not strongly connected. For each $i=2,...,d-2,$
the vertex $v_{i}$ is a critical one in $T(V)=T_{d-2,d-1}(v_{1},...,v_{d-1}).$
By the conditions $(a)$ and $(b)$ of item $(2),$ the vertex-sets
$\{v_{1},...,v_{i-1}\}$ and $\{v_{i+1},...,v_{d-1}\}$ cannot be joined
by a path without the use of $v_{i}$ in $T$ itself.
So, the conditions of Proposition 3 imply that
$V(T_{ncr})=\{w_{1},...,$
$w_{n-d+1}\}$ $=W$
and $V(T_{cr})=\{v_{1},...,v_{d-1}\}=V$
and hence,
the partition into two non-empty subsets $V=\{v_{1},...,v_{d-1}\}$ and
$W=\{w_{1},...,w_{n-d+1}\}$
with properties $(1)$ and $(2)$ is unique.

It is straightforward to check that $T_{d,n}$ is the tournament
with parameters $h_{1}=...=
h_{\lceil\frac{n-d+1}{2}\rceil}=d-2$
and
$h_{\lceil\frac{n-d+1}{2}\rceil+1}=...=h_{n-d+1}=1$ in ${\Cal H}_{d,n}.$
In turn, $T_{d,n}^{-}$ is the tournament
with parameters $h_{1}=...=
h_{\lfloor\frac{n-d+1}{2}\rfloor}=d-2$
and
$h_{\lfloor\frac{n-d+1}{2}\rfloor+1}=...=h_{n-d+1}=1$
in ${\Cal H}_{d,n}.$

In the sequel, we mainly consider the case $d=n-3.$
By Proposition 3, we have $V(T_{ncr})=W=
\{w_{1},w_{2},w_{3},w_{4}\}$ and
$T(W)=TT_{4}(w_{1},w_{2},w_{3},w_{4}).$
Obviously, $d_{T}^{+}(w_{i})=4-i+h_{i}$
for each $i=1,2,3,4.$ Since an isomorphism
takes non-critical vertices to non-critical vertices
and preserves their out-degrees and the sequence order of vertices
in the unique hamiltonian path of each transitive subtournament,
two elements $T$ and $T^{\prime}$
in the class ${\Cal H}_{n-3,n}$ are isomorphic if and only if
the out-degrees of $w_{2}$ and $w_{3}$ in $T$ (i.e. $h_{2}$ and $h_{3}$)
coincide with the out-degrees of the corresponding vertices
$w_{2}^{\prime}$ and $w_{3}^{\prime}$ in $T^{\prime}$
(i.e. $h_{2}^{\prime}$ and $h_{3}^{\prime}$). The parameters
$h_{2}$ and $h_{3}$ are related by $h_{3}\le h_{2}+1.$
By condition $(2a),$ for $h_{2}=n-5,$ it becomes $h_{3}\le h_{2}.$
Hence, denoting $h_{2}$ by $h$ for simplicity, we get that
the class
${\Cal H}_{n-3,n}$ consists of
$$\sum\limits_{h=1}^{n-6}(h+1)+(n-5)=\frac{n^{2}-7n+8}{2}$$
non-isomorphic tournaments.
Obviously, for any $T\in {\Cal H}_{n-3,n},$
we have $c_{n}(T)=1.$
According to Proposition 3 in [16],
if a tournament admits exactly one hamiltonian circuit,
then the same also holds for any of its
strongly connected subtournaments (as well as any subtournament
of $TT_{n}$ is also transitive and hence, contains exactly one hamiltonian path).
In particular, if $T\in {\Cal H}_{n-3,n},$ then $c_{\ell}(T)$
equals the number $s_{\ell}(T)$ of strongly connected subtournaments
of order $\ell$ in $T.$ Note that in [8], the statement
of Theorem II was formulated and proved precisely for the quantity
$s_{\ell}(T),$ but
for $c_{\ell}(T),$ it follows from this lower bound
on $s_{\ell}(T)$ and the Camion theorem. It is also easy to check
that the statement of Theorem III (as well as all the theorems
obtained in this paper) remains be true if one replaces $c_{\ell}(T)$
by $s_{\ell}(T)$ in its assumptions.

By Propositions 1 and 2,
if $T\in {\Cal T}_{\le n-3,n},$ where $n\ge 6,$ then
$|T_{ncr}|\ge 4$  and hence, $c_{n-1}(T)\ge 4.$
By definition, any $T\in {\Cal H}_{n-3,n}$ admits exactly $4$ non-critical vertices.
Hence, for each $T\in {\Cal H}_{n-3,n},$ we have
$c_{n-1}(T)=4.$ According to Theorem 2 [16], the converse is also true
for $T\in {\Cal T}_{\le n-3,n},$ where $n\ge 6.$

The case $\ell=n-2$ is much more complicated.
The point is that $c_{n-2}(T_{n-3,n})=\binom{4}{2}=6
> 5=n-(n-2)+3.$
However, exactly $6$ (not $5$) is
the minimum number of circuits of length $n-2$ in the class
${\Cal  T}_{\le n-3,n}.$
To show this, we need the following lemma on the relation between
$s_{n-2}(T)$ and $|T_{ncr}|.$

\smallskip

{\bf Lemma 3.} {\sl
Let $T$ be a strongly connected tournament of order $n\ge 9.$
Assume that there exists at most one non-critical vertex $w$
such that the corresponding strongly connected
subtournament $T-w$ is isomorphic to $T_{n-2,n-1}.$
Then $|T_{ncr}|\ge 4$ implies that $s_{n-2}(T)\ge 6.$
Moreover, if $|T_{ncr}|\ge 5,$ then $s_{n-2}(T)\ge 7.$}

{\sl Proof.} Assume that $T$
admits at least four non-critical vertices $w_{1},w_{2},w_{3},$ and $w_{4}.$
Set $W=\{w_{1},w_{2},w_{3},w_{4}\}.$
For each $i,$ denote by $W_{i}$ the set of non-critical
vertices in the strongly connected subtournament
$T-w_{i},$ i.e. the set of vertices $w$ for which
the subtournament $T-w_{i}-w$ is also strongly connected.
Obviously, if $w_{i}\in W_{j}$,
then $w_{j}\in W_{i}.$
By Corollary I, we have $|W_{i}|\ge 2$ for each $i.$
By the assumption of the lemma
and Theorem V, we can assume that if $|W_{i}|=2$ for some $i,$
then $i=1.$
Thus, $|W_{1}|\ge 2$ and
$|W_{i}|\ge 3$ for $i=2,3,4.$
Select a vertex $\tilde{w}$ in $W_{1}\setminus \{w_{4}\}.$
Obviously, without loss of generality we can assume that
$w_{3}\in W \setminus \{w_{1},w_{4},\tilde{w}\}.$
Then the strongly connected subtournaments $T-w_{4}-w,$
where $w\in W_{4},$ $T-w_{1}-\tilde{w},$ and
$T-w_{3}-\hat{w},$ where $\hat{w}\in W_{3}\setminus \{w_{4}\},$
are distinct and hence, as $|W_{3}|\ge 3$ and $|W_{4}|\ge 3,$
we have $s_{n-2}(T)\ge 6.$

Suppose now that there also exists a fifth non-critical vertex $w_{5}$ in $T.$
By the assumption, we have
$|W_{5}|\ge 3.$ Recall also that $|W_{2}|\ge 3.$
Obviously, at least one of the vertices $w_{2}$ and $w_{5}$
is different from $\tilde{w}.$
Hence, either $T-w_{2}-w,$
where $w\in W_{2}\setminus \{w_{3},w_{4}\},$ or $T-w_{5}-\hat{w},$
where $\hat{w}\in W_{5}\setminus \{w_{3},w_{4}\},$
is a $7$th strongly connected subtournament of order $n-2.$
In other words, $s_{n-2}(T)\ge 7$ if $|T_{ncr}|\ge 5.$
The lemma is proved. $\blacksquare$

\smallskip

Can one describe all strongly connected tournaments of order $n$
admitting at least two non-critical
vertices $w_{1}$ and $w_{2}$ for which the corresponding strongly
connected subtournaments $T-w_{1}$ and $T-w_{2}$ are isomorphic
to $T_{n-2,n-1}$? Obviously, removing the vertex $w_{1}$ or $w_{2}$ from
$T_{n-1,n}(w_{1},z_{1},...,z_{n-2},w_{2})$ leads us to
a strongly connected tournament of order $n-1$ and diameter $n-2.$
The same holds if we reverse the arc $(w_{2},w_{1}).$
Another example can be obtained from $T_{n-2,n-1}(z_{0},...,z_{n-2})$
by replacing the vertex $z_{i},$ where $i=0,...,n-2,$
with the transitive tournament $TT_{2}(\{w_{1},w_{2}\}).$

Note also that reversing the arcs $(z_{2},z_{3})$
and $(z_{3},z_{i}),$ where $i=0$ or $i=1,$ in
$T_{n-2,n-1}(z_{0},...,z_{n-2})$ yields
the tournament $T_{n-2,n-1}(z_{i+1},z_{i+2},z_{i},z_{3},...,z_{n-2}),$
where the subscript $i+2$ is considered with respect to
modulo $3.$ Hence, replacing the vertex $z_{3}$ by
$TT_{2}(\{w_{1},w_{2}\})$
in $T_{n-2,n-1}(z_{0},...,z_{n-2})$
and then reversing the arcs $(z_{2},w_{j})$ and $(w_{j},z_{i}),$
where $j=1$ and/or $j=2$ and
$i=0$ or $i=1$ (note that for different $j,$
we can take distinct possible values of $i$),
lead to a strongly connected tournament with two non-critical vertices
$w_{1}$ and $w_{2}$ such that $T-w_{1}\cong T-w_{2}\cong T_{n-2,n-1}.$
Of course, the same holds if we first replace $z_{n-5}$
by $TT_{2}(\{w_{1},w_{2}\})$
and then reverse the arcs $(w_{j},z_{n-4})$ and $(z_{i},w_{j})$
for $j=1$ and/or $j=2$ and
$i=n-3$ or $i=n-2$
(obviously, the above remark
on the possible choice of distinct $i$
for different $j$ remains true here).
One can show that
in fact, we have listed all the possible examples of strong
tournaments $T$ of order $n$ admitting at least two non-critical vertices
$w_{1}$ and $w_{2}$ for which $T-w_{1}$ and $T-w_{2}$
have diameter $n-2.$
However, in the sequel, we need only the following alternative
for the family of such tournaments.

\smallskip

{\bf Lemma 4.} {\sl
Let $T$ be a strongly connected tournament of order $n\ge 9.$
Assume that there exist two non-critical vertices $w_{1}$
and $w_{2}$ in $T$ such that the corresponding strongly connected
subtournaments $T-w_{1}$ and $T-w_{2}$ of order $n-1$
have the maximum possible diameter $n-2.$
Then either the diameter of the original tournament $T$ is not less than $n-2$
or $T$ admits at least $7$ strongly connected subtournaments of order $n-2.$}

{\sl Proof.} Obviously, we can assume that $T-w_{2}=T_{n-2,n-1}(z_{0},...,z_{n-2}).$
Then $w_{1}=z_{j}$ for some $j=0,...,n-2.$
Obviously, the cases $j=n-2,n-3,n-4,n-5$ are reduced to the cases
$j=0,1,2,3,$ respectively, by the transition
from the original tournament to its converse.
Hence, it suffices to consider only those $j$ for which
$0\le j\le n-6.$

Assume first that $j=0.$ Then $T-w_{2}-w_{1}=T_{n-3,n-2}(z_{1},...,z_{n-2}).$
Since $T-w_{1}\cong T_{n-2,n-1},$
either
$T-w_{1}=T_{n-2,n-1}(w_{2},z_{1},...,z_{n-2})$
or
$T-w_{1}=T_{n-2,n-1}(z_{1},...,z_{n-2},w_{2}).$
Recall that
$T-w_{2}=T_{n-2,n-1}(w_{1},z_{1},...,z_{n-2}).$
Hence, in the first case, for any orientation of the arc between
$w_{1}$ and $w_{2},$ the distance between $w_{2}$ and $z_{n-2}$ in $T$
is also equal to $n-2.$ In the second case, $T$ either coincides with
$T_{n-1,n}(w_{1},z_{1},...,z_{n-2},w_{2})$
(if $w_{2}\to w_{1}$) or is obtained from it by reversing
the arc from $w_{2}$ to $w_{1}$
(if $w_{1}\to w_{2}$). However, for the latter case,
there exist at least $7$ strongly connected subtournaments of order $n-2,$
namely, $T-w_{2}-w_{1},$ $T-w_{2}-z_{n-2},$ $T-w_{1}-z_{1},$
$T-z_{n-3}-z_{n-2},$ $T-z_{4}-z_{5},$ $T-z_{3}-z_{4},$
and $T-z_{1}-z_{2}.$

Assume now that $j=1.$
Then
$d_{T-w_{2}-w_{1}}^{+}(z_{0})=0$ and hence, $z_{0}\to w_{2},$
as $T-w_{1}$ is strongly connected. If $w_{2}\to z_{i}$ for some $i\ge 3,$
then the diameter of the subtournament $T-w_{1}$ is not greater than $n-3.$
Hence, $\{z_{3},...,z_{n-2}\}\Rightarrow w_{2}\to z_{2}.$
This means that
$T-w_{1}=T_{n-2,n-1}(z_{0},w_{2},z_{2},...,z_{n-2}).$
Since
$T-w_{2}=T_{n-2,n-1}(z_{0},w_{1},z_{2},...,z_{n-2}),$
for any orientation of the arc between
$w_{1}$ and $w_{2},$ the distance between $z_{0}$ and $z_{n-2}$ in $T$
is also equal to $n-2.$

Finally, assume that $2\le j\le n-6.$
Obviously, for this case,
$d_{T-w_{2}-w_{1}}^{+}(z_{j+1})=$ $d_{T-w_{2}-w_{1}}^{+}(z_{j+2})=j+1.$
If $\{z_{j+1},z_{j+2}\}\Rightarrow w_{2}$
or  $\{z_{j+1},z_{j+2}\}\Leftarrow w_{2},$
then we have $d_{T-w_{1}}^{+}(z_{j+1})=d_{T-w_{1}}^{+}(z_{j+2})=j+2$
or $d_{T-w_{1}}^{+}(z_{j+1})=d_{T-w_{1}}^{+}(z_{j+2})=j+1,$
respectively.
But this is impossible in
$T_{n-2,n-1}$ because $2\le j\le n-6.$
In turn, if $z_{j+1}\to w_{2}\to z_{j+2},$
then $d_{T-w_{1}}^{+}(z_{j+1})=j+2$ and $d_{T-w_{1}}^{+}(z_{j+2})=j+1.$
Since $T-w_{1}\cong T_{n-2,n-1},$ we must have
$z_{j+2}\to z_{j+1},$ a contradiction.
Hence, $z_{j+2}\to w_{2}\to z_{j+1}.$

Since $T-w_{2}=T_{n-2,n-1}(z_{0},z_{1},...,z_{n-3},z_{n-2}),$
both $T-w_{2}-z_{0}$ and
$T-w_{2}-z_{n-2}$ are strongly connected.
Moreover, removing any of $2$-element vertex-sets
$\{z_{0},z_{1}\},$ $\{z_{n-3},z_{n-2}\},$
and $\{z_{0},z_{n-2}\}$ from $T-w_{2}$
yields a strongly connected tournament of order $n-3.$
Since $z_{j+2}\to w_{2}\to z_{j+1},$ where $2\le j\le n-6,$
this implies that
$T-z_{0}-z_{1},$
$T-z_{n-3}-z_{n-2},$ and
$T-z_{0}-z_{n-2}$
are also strongly connected.
Obviously, by the condition $w_{1}=z_{j},$ where $2\le j\le n-6,$
they are all different from two strongly connected subtournaments
of the form $T-w_{1}-w,$ where $w$ is one of the two
non-critical vertices of the strongly connected tournament $T-w_{1}.$
Thus, for the considered case, we have $s_{n-2}(T)\ge 7.$
The lemma is proved. $\blacksquare$

\smallskip

Now we are in a position to formulate and prove the main result of this section
which implies that in fact, the statement of Conjecture 1 also holds
for $d=n-3$ and $\ell=n-2.$

\smallskip

{\bf Theorem 4.} {\sl
Let $T$ be a strongly connected tournament of order $n\ge 9$ whose diameter
is not greater than $n-3.$
Then $c_{n-2}(T)\ge c_{n-2}(T_{n-3,n})=6$ with equality holding
if and only if $T\cong T_{n-3,n}.$}

{\sl Proof.} Indeed, by Propositions 1 and 2,
any tournament $T$ in the class ${\Cal T}_{\le n-3,n},$
where $n\ge 9,$
contains at least four non-critical vertices.
Assume that $c_{n-2}(T)\le 6.$
Then, by Lemma 4,
there exists at most one non-critical vertex $w$
for which $T-w\cong T_{n-2,n-1}$ and hence,
Lemma 3 implies that $c_{n-2}(T)=6$ and
$|T_{ncr}|=4.$
In particular, according to Propositions 1 and 2, again, we have
$T\in \Cal{T}_{n-3,n},$ where $n\ge 9.$
If the subtournament $T_{ncr}$ is strongly connected, then,
by Lemma 1, either
$T=T_{ncr}$ or
$T=\Delta(v_{1},v_{2},T_{ncr}).$
This implies that either $|T|=|T_{ncr}|=4 <9$
or $|T|=|T_{ncr}|+2=6<9.$
Hence, in the sequel, we can assume that the subtournament
$T_{ncr}$ is not strongly connected.

Assume first that the subtournament $T_{ncr}$
contains exactly two strongly connected components $T_{1}$ and $T_{2}.$
Without loss of generality, we can assume that $T_{1}\Rightarrow T_{2}.$
By Corollary 1, we have
$T=T_{n-3,n-2}(T_{2},v_{1},...,v_{n-4},T_{1}).$
Since there exists no strongly connected tournament
of order $2,$ either
$|T_{1}|=3$ and $|T_{2}|=1$
or
$|T_{1}|=1$ and $|T_{2}|=3.$
By duality, it suffices to consider only the first case.
Assume that $V(T_{2})=\{w\}.$
By Lemma 2, for any two non-critical vertices $w_{1}$ and $w_{2},$
the subtournament $T-w_{1}-w_{2}$ is also strongly connected.
Moreover, $T-w-v_{1}=T_{n-5,n-4}(v_{2},...,v_{n-4},T_{1}).$
Hence, for the considered case, the inequality $c_{n-2}(T)\ge \binom{4}{2}+1
=7$ holds.

Finally, if $T_{ncr}$
contains at least three strongly connected components, then
$T_{ncr}=TT_{4}(w_{1},$ $w_{2},w_{3},w_{4}).$
Let
$Q=w_{4},v_{1},...,v_{t},w_{1}$ be a shortest path from $w_{4}$ to $w_{1}$
in $T.$
Then the vertices of the closed walk $w_{4},v_{1},...,v_{t},
w_{1},w_{2},w_{3},w_{4}$ induce a strongly connected subtournament $T_{Q}$
containing all the non-critical vertices. Since $T_{ncr}$ is not strongly connected,
Lemma 1 implies that $T_{Q}=T$ and hence,
the path $Q$ must include all $n-4$ critical vertices of $T.$ This means that $t\ge n-4.$
On the other hand, $t\le n-4$ because the distance between $w_{4}$ and $w_{1}$
cannot be greater than $n-3.$
Thus, $t=n-4$ and
$\{v_{1},...,v_{n-4}\}$ is the set of critical vertices
of the strongly connected tournament $T.$

Since
$w_{4},v_{1},...,v_{n-4},w_{1}$ is a shortest path from $w_{4}$ to $w_{1}$ in $T,$
we have $\{v_{2},...,$ $v_{n-4}\}$ $\Rightarrow w_{4}\to v_{1}$
and
$v_{n-4}\to w_{1}\Rightarrow\{v_{1},...,v_{n-5}\}.$
Suppose that $w_{3}\Rightarrow\{v_{1},...,$ $v_{n-4}\}.$
If $w_{2}\Rightarrow\{v_{1},...,v_{n-5}\},$
then a shortest path from $w_{4}$ to $w_{3}$ in $T$
has the form $w_{4},v_{1},...,v_{n-4},w_{1},w_{3}$ or
$w_{4},v_{1},...,v_{n-4},w_{2},w_{3}$
and hence, the distance between $w_{4}$ and $w_{3}$ in $T$ equals $n-2,$
which is impossible.
This implies that $v_{p}\to w_{2}$ for some $1\le p\le n-5.$
If $p\le n-6,$ then the cycle
$w_{4},v_{1},...,v_{p},w_{2},w_{3},v_{n-4},w_{1},w_{4}$
contains all the non-critical vertices but does not include the critical vertex $v_{n-5},$ which is impossible
in the case where $T_{ncr}$ is not strongly connected by Lemma 1.
Hence, $v_{n-5}\to w_{2}\Rightarrow\{v_{1},...,v_{n-6}\}.$ However, for this case,
$T-w_{4}-v_{1}$ ($v_{2},...,v_{n-4},w_{1},w_{2},w_{3},v_{2}$)
and $T-w_{1}-v_{n-4}$ ($v_{1},...,v_{n-5},w_{2},w_{3},w_{4},v_{1}$)
together with
$T-w_{4}-w_{3}$  ($v_{1},...,v_{n-4},w_{1},w_{2},v_{1}$),
$T-w_{4}-w_{2}$  ($v_{1},...,v_{n-4},w_{1},w_{3},v_{1}$),
$T-w_{4}-w_{1}$  ($v_{1},...,v_{n-4},v_{n-6},v_{n-5},w_{2},w_{3},v_{1}$),
$T-w_{3}-w_{2}$  ($v_{1},...,v_{n-4},w_{1},w_{4},v_{1}$),
and
$T-w_{3}-w_{1}$ ($v_{1},...,v_{n-4},v_{n-6},v_{n-5},w_{2},w_{4},$ $v_{1}$)
are strongly connected tournaments of order $n-2$ (with the spanning closed
walks given in the parentheses).
This shows that the case $w_{3}\Rightarrow\{v_{1},...,v_{n-4}\}$ is impossible.
Obviously, if $\{v_{1},...,v_{n-4}\}\Rightarrow w_{3},$
then the subtournament $T-w_{4}$ is not strongly connected. Hence,
$v_{p}\to w_{3}\to v_{k}$
for some $1\le p,k\le n-4.$
Similarly,
$v_{m}\to w_{2}\to v_{s}$
for some $1\le m,s \le n-4.$

As $V(T_{cr})=\{v_{1},...,v_{n-4}\}$ and
$v_{1},...,v_{n-4}$ is a shortest path from $v_{1}$ to $v_{n-4},$
we have
$T_{cr}
=T_{n-5,n-4}(v_{1},...,v_{n-4}).$ In particular, it is strongly connected.
Since, as we have seen above,
for any $w\in V(T_{ncr}),$ there exist both an arc from $T_{cr}$ to $w$
and an arc from $w$ to $T_{cr},$
adding any two of the four non-critical
vertices to $T_{cr}$ yields a strongly connected tournament of order $n-2.$
If $w_{3}\to v_{k}$ for some $k\ge 2,$
then $w_{3},v_{k},...,v_{n-4},v_{2},...,v_{n-4},w_{1},w_{2},w_{3}$
is a spanning closed walk in $T-w_{4}-v_{1}.$ This means that
$T-w_{4}-v_{1}$ is a $7$th strongly connected subtournament of order $n-2$ in $T,$
which is impossible.
Hence, $\{v_{2},...,v_{n-4}\}\Rightarrow w_{3}\to v_{1}.$
Similarly, $v_{n-4}\to w_{2}\Rightarrow\{v_{1},...,v_{n-5}\}.$
This implies that $T=T_{n-3,n}(\widehat{w_{3},w_{4},v_{1}},...,\widehat{v_{n-4},w_{1},w_{2}}).$
The theorem is proved. $\blacksquare$

\smallskip

Theorem 4 implies that the following conjecture first presented in [16]
is true at least for $d=n-3$ and $h=2.$

\smallskip

{\bf Conjecture 3 [16]}.
{\sl Let $T$ be a strongly connected tournament of order $n\ge 4$

with diameter $d\ge 3.$ If $2h\le n-d+1,$ then}
$$c_{n-h}(T)\ge \binom{n-d+1}{h}.$$

{\sl Moreover, for sufficiently large $n,$
the equality holds if and only if the following

two conditions are satisfied:

(1) there exist two
vertex-sets $V=\{v_{1},...,v_{d-1}\}$ and $W=\{w_{1},...,w_{n-d+1}\}$

such that

(a) $T(V)=T_{d-2,d-1}(v_{1},...,v_{d-1}),$

(b) $T(W)=TT_{n-d+1}(w_{1},...,w_{n-d+1});$

(2) for each $i=1,...,n-d+1,$ there exists $h_{i}$
such that
$\{v_{h_{i}+1},...,v_{d-1}\}\Rightarrow w_{i}\Rightarrow$

$\{v_{1},...,v_{h_{i}}\}$ with

(a) $h_{1}=...=h_{h}=d-2,$
$1\le h_{i}\le d-2$ for $i=h+1,...,n-d+1-h,$

and $h_{n-d+2-h}=...=h_{n-d+1}=1,$

(b) $h_{s}\le h_{r}+1\ \text{ for any }\ h+1\le r\le s \le n-d+1-h.$}

\smallskip

{\bf Remark 4}. If $h\le n-d$ and $T\in {\Cal H}_{d,n},$
then removing any $h$ non-critical vertices from $T$
yields a strongly connected subtournament of order $n-h.$
Assume now that $2h\le n-d+1$ and $T\in {\Cal H}_{d,n}.$
Then any strongly connected subtournament of order $n-h$ in $T$
has this form
if and only if $T$ satisfies condition $(2a)$
(all the other conditions of Conjecture 3
are automatically satisfied for $T\in {\Cal H}_{d,n}$).
Obviously, for this case, $s_{n-h}(T)$ (as well as $c_{n-h}(T)$)
equals the combination of $n-d+1$ things $h$ at a time.
Hence, we can say that
the statement of Conjecture 3 holds for the class ${\Cal H}_{d,n}.$
The case of arbitrary $T\in {\Cal T}_{d,n}$  is much more complicated
and, as a consequence, is still open.

\smallskip

In conclusion, we consider the case of small values of the length $\ell\ge 3.$
The main result of [3] implies that
if $T\in {\Cal H}_{n-3,n},$ where $n\ge 6,$ and $h_{3}\le h_{2},$
then $c_{3}(T)=n-3+1=n-2 <n=n-3+3.$
Indeed, for each $i=2,3,$ only
$w_{i},v_{h_{i}},v_{h_{i}+1},w_{i}$
is a circuit of length $3$ containing the vertex $w_{i}$
(note that if $h_{3}\ge h_{2}+1,$ then
$w_{3},v_{h_{3}},w_{2},w_{3}$
is also a circuit of length $3$ passing through $w_{3}$)
and hence, $c_{3}(T)=
c_{3}(T-w_{2}-w_{3})+c_{3}(T,w_{2})+c_{3}(T-w_{2},w_{3})=
c_{3}(T_{n-3,n-2})+2*1=(n-2)-3+1+2=n-2.$
Thus, in ${\Cal T}_{n-3,n},$
there are at least (in fact, according to [3], exactly)
$$\sum\limits_{h=1}^{n-5}h=\frac{(n-4)(n-5)}{2}$$
non-isomorphic tournaments having precisely $c_{3}(T_{n-3,n})$
circuits of length $3.$

In turn,
if $T\in {\Cal H}_{n-3,n}$ and
$h_{3}\le h_{2}-1$ (this is possible iff $n\ge 7$),
then for each $i=2,3,$ only
$w_{i},v_{h_{i}-1},v_{h_{i}},v_{h_{i}+1},w_{i}$ and
$w_{i},v_{h_{i}},v_{h_{i}+1},v_{h_{i}+2},w_{i}$ (here we assume that
$v_{0}=w_{4}$ and $v_{n-3}=w_{1}$)
are circuits of length $4$ containing $w_{i}$
(note that if $h_{3}\ge h_{2},$ then
$w_{3},v_{h_{3}},v_{h_{3}+1},w_{2},w_{3}$
is also a circuit of length $4$ passing through $w_{3}$)
and hence,
$c_{4}(T)=
c_{4}(T-w_{2}-w_{3})+c_{4}(T,w_{2})+c_{4}(T-w_{2},w_{3})=
c_{4}(T_{n-3,n-2})+2*2=
(n-2)-4+1+4=n-1=n-4+3.$
Thus, for $n\ge 7,$ in ${\Cal T}_{n-3,n},$
there are at least (according to our conjecture, precisely)
$$\sum\limits_{h=1}^{n-5}(h-1)=\frac{(n-6)(n-5)}{2}$$
non-isomorphic tournaments having exactly $c_{4}(T_{n-3,n})$
circuits of length $4.$
In the previous section, just after the proof of Theorem 2,
we have showed that the inequality $c_{4}(T)\le n-2$ implies that
the diameter of $T$ is not less than $n-2.$
Thus, as we have stated in Conjecture 1,
the minimum of
$c_{4}(T)$ in the class $\Cal{T}_{\le n-3,n}$ is attained at
$T_{n-3,n}$ and equals $n-1,$ indeed.

Finally, let us consider the case $\ell=5.$ As it was shown in the
very beginning of the proof of Theorem 3,
for $T\in {\Cal T}_{\le n-3,n},$ we always have
$c_{5}(T)\ge n-5+3=n-2.$
Obviously, if $T=T_{n-3,n}(\widehat{w_{3},w_{4},v_{1}},...,\widehat{v_{n-4},w_{1},w_{2}}),$
then $c_{5}(T-w_{3})=c_{5}(T_{n-3,n-1})=(n-1)-5+2=n-4.$
Moreover, $w_{3},v_{1},v_{2},v_{3},v_{4},w_{3}$ and
$w_{3},w_{4},v_{1},v_{2},v_{3},w_{3}$ are the only cycles
of length $5$ passing through the vertex $w_{3}.$
Thus, $c_{5}(T)=c_{5}(T-w_{3})+c_{5}(T,w_{3})=n-4+2=n-2.$

Repeating practically word for word the proof of Theorem 3 for the case $\ell=5,$
one can show that the equality $c_{5}(T)=n-2$ in the class ${\Cal T}_{\le n-3,n}$
implies that $T\cong T_{n-3,n}.$
Indeed, we need to comment only the case
$z_{k}\to w\to z_{k+1}$ for $3\le k\le n-6$ because the cycle
$\gamma_{i}=z_{k-\ell+2+i},...,z_{k},w,z_{k+1},...,z_{k+i},z_{k-\ell+2+i}$
makes no sense for $\ell=5$ and $i=4.$
However, it can be constructed in the same way as it was done
for the cycle $\gamma_{3}$ of length $4$ just after the proof of
Theorem 2 if one
considers
the subtournament $S_{5}$ induced by the vertex-set
$V_{5}=\{z_{k-3},z_{k-2},z_{k-1},z_{k},w\}$ of order $5.$

The fact that the statement of Theorem 3 is also true for $\ell=5$ and
for $\ell=n-2$
(in the latter case, after replacing $n-\ell+3$ by $n-\ell+4$)
implies that the statement of Conjecture 1 can be strengthened.
This will be done in one of our further papers.

\smallskip

{\bf Acknowledgements}

The author is grateful to the Editorial Board
(personally to A.A. Schkalikov)
for the fact that it was able to find a way out of
the situation that had developed with
the paper by the spring of 2015.
He is also thankful to the referee for his careful reading of the work
and his valuable comments
which we hope greatly improved
the original text
and thereby made it more understandable for the reader.

\bigskip

\centerline{\bf References }

\smallskip

[1] J. Bang-Jensen  and  G. Gutin,
{\sl Digraphs: Theory, Algorithms and Applica- tions},
Springer-Verlag, London, 2000.

[2] M. Burzio and D.C. Demaria,
On a classification of Hamiltonian tournaments,
{\sl Acta Univ. Carolin. - Math. Phys.,} Prague, {\bf 29} (2) (1988), 3-14.

[3] M. Burzio and D.C. Demaria,
Hamiltonian tournaments with the least number of $3$-cycles,
{\sl J. Graph Theory} {\bf 14} (1990), 663-672.

[4] A.H. Busch, A note on the number of hamiltonian paths in
strong tournaments, {\sl Electron. J. Combin.} {\bf 13} (2006), $\#N3.$

[5] P. Camion, Chemins et circuits hamiltoniens des graphes complets,
{\sl C.R. Acad. Sci. Paris } {\bf 249} (1959), 2151-2152.

[6] R.J. Douglas, Tournaments that admit exactly one hamiltonian circuit,
{\sl Proc. London Math. Soc.} {\bf 21} (1970), 716-730.

[7] M. Las Vergnas,  Sur le nombre de circuits dans un
tournoi fortement connexe,
{\sl Cahiers Centre \'Etudes Recherche Op\'er.} {\bf 17} (1975), 261-265.

[8] J.W. Moon,  On subtournaments of a tournament,
{\sl Canad. Math. Bull.} {\bf 9} (1966), 297-301.

[9] J.W. Moon, {\sl Topics on Tournaments},
Holt, Rinehart and Winston, New York, 1968.

[10] J.W. Moon, The minimum number of spanning paths in a strong tournament,
{\sl Publ. Math. Debrecen} {\bf 19} (1972), 101-104.

[11] G.V. Nenashev, On the existense of non-critical vertices in digraphs,
{\sl Zap. POMI} {\bf 406} (2012), 107-116.
English trans.: J. Math. Sci. (N.Y.) {\bf 196} (6) (2014),
791-796.

[12] S.B. Rao and A.R. Rao, The number of cut vertices
and cut arcs in a strong directed graph,
{\sl Acta Math. Sci. Hungar.} {\bf 22} (1971), 411-421.

[13] L. R\'edei, Ein kombinatorischer Satz,
{\sl Acta Litt. Sci. Szeged} {\bf 7} (1934), 39-43.

[14] S.V. Savchenko, On the number of non-critical vertices
in strongly connected digraphs, {\sl Mat. Zametki} {\bf 79} (2006), 743-755.
English trans.: Math. Notes {\bf 79} (2006), 687-696.

[15] S.V. Savchenko, On the number of non-critical vertices
in strong tournaments of order $n$
with minimum out-degree $\delta^{+}$ and in-degree $\delta^{-},$
{\sl Discrete Math.} {\bf 310} (2010), 1177-1183.

[16] S.V. Savchenko, Non-critical vertices and long circuits
in strong tournaments of order $n$ and diameter $d,$
{\sl J. Graph Theory} {\bf 70} (2012), 361-383.

[17] C. Thomassen,  Whitney's 2-switching theorem,
cycle spaces, and arc mappings of directed graphs,
{\sl J. Combinat. Theory Ser. B} {\bf 46} (1989), 257-291.

[18] C. Thomassen, On the number of Hamiltonian cycles in tournaments,
{\sl Discrete Math.} {\bf 31} (1980), 315-323.

\end{document}